\documentclass[12pt,a4paper]{article}
\usepackage[latin1]{inputenc}
\usepackage{amsmath,amssymb}

\newcommand{\mltsum}[2]{\sum\limits_{\stackrel{\scriptstyle #1}{#2}}}

\author{S.J.Patterson}
\title{Selberg Sums -- a new perspective}
\date{}

\begin{document}
\maketitle

\begin{abstract}  
Selberg sums are the analogues over finite fields of certain integrals studied
by Selberg in in 1940s.   The original versions of these sums were introduced by
R.J.Evans in 1981 and, following an elegant idea of G.W.Anderson in 1991 they
were evaluated by Anderson, Evans and P.B.~van~Wamelen.  In 2007 the author 
noted that these sums and certain generalizations of them appear in the study 
of the distribution of Gauss sums over a rational function field over a finite
field.   The distribution of Gauss sums is closely related to the distribution 
of the values of the discriminant of polynomials of a fixed degree.  Here we 
shall take this up further.   The main goal here is to establish the basic 
properties of Selberg sums and to formulate the problems which arise from this
point of view.
\end{abstract}
\section{Introduction}
The usual class of character sums known as Selberg sums was introduced in 1981 by
R.J.Evans, \cite{RJE1}.   They were the analogues over finite fields of a class,
or rther several related classes, of integrals introduced by A. Selberg in 1944, 
\cite[p.204 ff.]{ASCW1}. In fact he had already used an integral of this type in 
1941 in \cite[p.74 ff.]{ASCW1} but was uncertain as to whether the integral was 
already known and waited before publishing details.   

Characteristic of these integrals is that they are taken over a space of monic 
polynomials of a fixed degree and that one factor in the integrand is a power of 
the discriminant.   Selberg regarded the integrals he studied as a extension of 
Euler's beta-function and, indeed, his evaluation gives the integrals as a quotient 
of products of gamma functions.  The primary example is the integral
\[
\int_0^1\cdots \int_0^1 (x_1\cdots x_i)^{\alpha-1}
((1-x_1)\cdots (1-x_n))^{\beta-1} |\prod_{i < j} (x_i-x_j)|^{2\gamma}
\mathrm dx_1 \cdots \mathrm dx_i
\]
in the region $\mathrm{Re}(\alpha)>0$, $\mathrm{Re}(\beta)>0$,
$\mathrm{Re}(\gamma) > -\mathrm{Min}(\frac 1i, \frac{\mathrm{Re}(\alpha)}{i-1},
\frac{\mathrm{Re}(\beta)}{i-1})$.
Selberg's evaluation shows that the integral is equal to 
\[
\prod\limits_{j=1}^{i} \frac{\Gamma(1+j\gamma)\Gamma(\alpha+(j-1)\gamma)\Gamma(\beta+(j-1)\gamma)}
{\Gamma(1+\gamma)\Gamma(\alpha+\beta+(i+j-2)\gamma)}.
\]
For our purposes a transformation is helpful.    This integral can be regarded as over all 
monic polynomials of degree $i$ with all their roots in $[0,1]$.   We let $\sigma_1,\dots,\sigma_i$ 
be the standard symmetric functions in $x_1,\dots,x_i$.   Then one has 
\[
\mathrm d\sigma_1\wedge \mathrm d\sigma_2\wedge
\dots\wedge\mathrm d\sigma_i =
\prod_{j<j'} (x_j-x_{j'})\ 
\mathrm d x_1\wedge \mathrm dx_2 \wedge \dots\wedge
\mathrm dx_i.
\]
It follows that we can write Selberg's integral as 
\[
\int |f(0)|^\alpha |f(1)|^\beta |D(f)|^{\gamma-\frac 12} 
\mathrm df,
\]
where the integral is over the set of polynomials described above and $\mathrm df$
is the standard Lebesgue measure on the affine space of the $f$.  

We can now describe Evans' analogue.   Let $q$ be a power of the odd prime number $p$.
Let $\chi_1,\chi_2$ and $\chi_3$ be multiplicative characters on $\mathbb F_q^\times$ 
and let $\omega$ be the unique quadratic character.   Then Evans introduces 
\[
\sum \chi_1(f(0))\chi_2(f(1))(\omega \chi_3)(D(f))
\] 
where $D(f)$ denotes the discriminant of the polynomial $f$ and the sum is taken over 
all monic polynomials of a fixed degree $i$.    It turns out 
that introducing the factor $\omega$ makes the development much smoother and it 
can be regarded as the analogue of the term $-\frac 12$ in the exponent of 
$D(f)$ in the alternative version of Selberg's integral.  

In 1990 G.W. Anderson \cite{GWA} found an ingenious and remarkably simple method for evaluating 
some of these sums.   Following this first R.J.Evans, \cite{RJE2} and then P.B. v. Wamelen \cite{vW}
extended Anderson's method to cover all cases.   Details of the results are given in \S 3.

Since Selberg's integral can be considered as an extension of the beta-integral we can regard 
Evans' sum as an extension of their analogues over finite fields, namely Jacobi sums.  One 
immediate generalization of the beta-integral is the standard integral representation of the 
hypergeometric function and their generalizations.   In \cite{P1} it was pointed out that 
the finite field analogues of these functions intervene in the theory of metaplectic forms,
and especially of Eisenstein series, over rational function fields over finite fields.   In
this article we shall follow this line of thought and apply the appelation ``Selberg sum'' to 
this more general class of sum.   To the best of my knowledge the corresponding archimedean 
functions have not been investigated.

The theory of metaplectic forms leads us to the evaluation to some more general Selberg sums.
This is an aspect which we shall not go into here; some indications are given in \cite{P2}.
He we shall take up a theme of that paper, namely a transformation formula.   This will be given
in Section 4.   It is quite elementary but among other things it shows that the results of 
Anderson, Evans and v.Wamelen can be used to evaluate further interesting classes of Selberg
sums.

It is worth noting that it is relatively easy to compute specific examples of Selberg sums and
so one can investigate them experimentally.   We shall explain why this is so later.   This means 
that one can investigate the sums experimentally and, for example, examine their sizes in various
metrics.   This is a topic which we shall postpone to a later paper.

My thanks are due to Ben Brubaker, Paul Gunnells and Zeev Rudnick for helpful conversations 
about the topics discussed in this paper.

\section{Notations}

Let, as before, $q$ be a power of the odd prime $q$.   For $f,g \in \mathbb F_q[x]$ let $R(f,g)$ 
be the resultant of $R(f,g)$  We recall that $R(f,g)$ is a universal polynomial (Sylvester's 
determinant) in the coefficients of $f$ and $g$.   The function $g \mapsto R(f,g)$ depends only
on $g \pmod f$.   Also $R(f,g)$ is bimultiplicative and satisfies the reciprocity law
\[
R(g,f)=(-1)^{\mathrm{deg}(f)\mathrm{deg}(g)}R(f,g).
\] 
Finally $R(x-a,f)=f(a)$.   These properties determine $R$ completely and can be cast into 
the form of an efficient algorithm based on the model of continued fractions for its evaluation.

For $f \in \mathbb F_q[x]$ let $D(f)$ be the discriminant of $f$.   If $f$ is monic one has 
$D(f)=\eta(\mathrm{deg}(f)) R(f,f')$ where $\eta(i) = 1$ for $i \equiv 0,1 \pmod 4$ and 
$\eta(i) = -1$ for $i \equiv 2,3 \pmod 4$\footnote{Note that this factor is missing from the 
formula (5) on \cite[p.87]{vdW}, an oversight which unfortunately led to author to similar 
inaccuracies in \cite{P1}.} .   

We let $\mu$ be the M\"obius function on  $\mathbb F_q[x]$.   Then Pellet's formula (see 
\cite[(1.5)]{CR}, \cite[Lemma 4.1]{KC} or \cite{RGS}) asserts that
\[
\mu(f) = (-1)^{\mathrm{deg}(f)}\omega(D(f))
\]
where, as above, $\omega$ is the unique quadratic character on $\mathbb F_q^\times$.

Let $X_q$ be the group of characters on $\mathbb F_q^\times$ which we consider as 
taking values in $\mathbb Q(1^{\frac 1 {q-1}})$.    For $\chi \in X_q$ let $\mathrm{ord}(\chi)$
denote the order of $\chi$.   Let $e_o$ be a non--trivial additive character on $\mathbb F_q$;
for our purposes it is convenient to regard it as taking values in $\mathbb Q(1^{\frac 1 p})$.
The choice of $e_o$ will not play any significant role in our discussion.   For $\chi \in X_q$
we define the Gauss sum over $\mathbb F_q$ to be 
\[
\tau(\chi) = \sum_{a \in \mathbb F_q^\times}\chi(a)e_o(a).
\]
For $\chi_1,\chi_2 \in X_q$ let 
\[
J(\chi_1,\chi_2) = \mltsum{a \in \mathbb F_q}{a\neq 0,1}  \chi_1(a)\chi_2(1-a)
\]
be the Jacobi sum.   Then if $\chi_1 \chi_2 \neq 1$ 
\[
\tau(\chi_1)\tau(\chi_2) = J(\chi_1,\chi_2)\tau(\chi_1\chi_2)
\]
whereas if  $\chi_1 \chi_2 = 1$, but $\chi_1, \chi_2 \neq 1$ then
\[
\tau(\chi_1)\tau(\chi_2) = \chi_1(-1)q
\]
and 
\[
J(\chi_1,\chi_2) = -\chi_1(-1).
\]
Additionally one has $J(\chi,1)=-1$ if $\chi \neq 1$ and  $J(1,1)=q-2$.   Also 
$J(\chi_1,\chi_2)=J(\chi_2,\chi_1)$, $J(\chi_1,\chi_2)=\chi_2(-1)J((\chi_1\chi_2)^{-1},\chi_2)$
and $J(\chi_1,\chi_2)=\chi_1(-1)J(\chi_1,(\chi_1\chi_2)^{-1})$.  There are some further 
relations and properties of these sums which have been studied in great detail.

Let now $k=\mathbb F_q(x)$ and $R=\mathbb F_q[x]$.   We define 
\[
e:k \rightarrow \mathbb Q(1^{\frac 1 p})
\]
by 
\[
\begin{array}{rl}
e(f) &=  e_o(\sum_{v \neq \infty} \mathrm{Res}_v (f \mathrm dx))\\
     &=  e_o(-\mathrm{Res}_\infty (f \mathrm dx)).
\end{array}
\]
Here, in the sum, $v$ runs through the finite places of $k$.

We define for $\chi \in X_q$ and $g \in R-\{0\}$ the Dirichlet character 
of modulus $g$ 
\[
f \mapsto \chi(f/g) = \chi(R(g,f)).
\]
We define the global Gauss sum associated with $k$ to be
\[
g(r,\chi,c) = \sum_{d \pmod c} \chi(d/c) e(rd/c)
\]
for $r \in R-\{0\}$.   This can, by means of the Davenport-Hasse theorem,
be evaluated in terms of Gauss sums over $\mathbb F_q$.   If $r$ is coprime 
to $c$ one obtains
\[
g(r,\chi,c)= \mu(c) \chi(r/c)^{-1}\chi(c'/c)(-\tau(\chi))^{\mathrm{deg}(c)}.
\]
Here $c'$ denotes the derivative of $c$.   This shows that the function $\chi(c'/c)$
is as subtle a function of $c$ as $g(r,\chi,c)$.   Recall that it is equal to 
$\chi(\eta(\mathrm{deg}(c)))\chi(D(c))$ for $c$ monic.  We therefore have for $c$ coprime
to $r$
\[
g(r,\chi,c)= \mu(c) \chi(\eta(\mathrm{deg}(c)))\chi(r/c)^{-1}(-\tau(\chi))^{\mathrm{deg}(c)}\chi(D(c)).
\]
In the special case $\chi=\omega$, $r=1$ we obtain
\[
g(r,\omega,c)= \mu(c)^2 \omega(\eta(\mathrm{deg}(c)))\tau(\omega)^{\mathrm{deg}(c)}.
\]
As $\tau(\omega)$ is $\sqrt q$ times a fourth root of $1$ this gives us the Gauss evaluation of the 
quadratic Gauss sum in this context.   One could also use this argument to prove Pellet's formula 
using that evaluation, as given, for example, in \cite[XIII,\S 11, Proof of Theorem 13]{BNT}.

We can now define the Selberg sums which we are going to investigate here.   Apart from a factor 
$\pm 1$ it specializes to the one used in the papers of Anderson, Evans and v.~Wamelen when we 
take $r$ to be of the form $x^{e_0}(x-1)^{e_1}$.   For $\chi_1,\chi_2$ it is 
\[
\sum_c \chi_1(r/c)\omega(D(c)) \chi_2(D(c))
\]
where $c$ is summed over all monic polynomials of a fixed degree $i$.   For our purposes it is
convenient to modify the expression by using Pellet's formula.   It becomes
\[
(-1)^i\sum_c \mu(c) \chi_1(r/c) \chi_2(D(c)),
\]
or, with the reciprocity law,
\[
(-1)^i\chi(-1)^{i\mathrm{deg}(r)i}\sum_c \mu(c) \chi_1(c/r) \chi_2(D(c)).
\]
Both of the sums here could be considered as to be entitled to be the more fundamental 
one.   We choose
\[
Se(r,\chi_1,\chi_2,i)= \mltsum{c\mbox{ \scriptsize{monic}}}{\mathrm{deg}(c)=i} \mu(c) \chi_1(r/c) \chi_2(D(c)).
\]
We take $\chi_1(r/c)$  to be zero if $r$ and $c$ are not coprime.   By means of the 
Davenport-Hasse theorem we have then 
\[
Se(r,\chi_1,\chi_2,i)= (-1)^i \tau(\chi_2)^{-i} \chi_2(\eta(i)) \mltsum{c\mbox{ \scriptsize{monic}}}{\mathrm{deg}(c)=i} \chi_1(r/c) 
g(1,\chi_2,c).
\]
If now there exists an exponent $e$ so that 
\[
\chi_1\chi_2^e =1
\] 
then 
\[
Se(r,\chi_1,\chi_2,i)= (-1)^i \tau(\chi_2)^{-i} \chi_2(\eta(i)) \mltsum{c\mbox{ \scriptsize{monic}}}{\mathrm{deg}(c)=i} 
g(r^e,\chi_2,c)
\]
where we restrict the sum to $c$ coprime to $r$.   In this case one can use of the theory of metaplectic
forms to investigate this sums, see \cite{P1},\cite{P2}.   At the other extreme $\chi_2=1$ and 
$Se(r,\chi_1,1,i)$ is then the coefficient of $q^{-is}$ in the L-series $L(s,\chi_1(r/\cdot))$ (over $R$). 
It follows from this that if $\chi_1(r/\cdot)$ is non-principal then $Se(r,\chi_1,1,i)=0$ for 
$i\ge \mathrm{deg}(r_0)$ where $r_0$ is the conductor of $\chi_1(r/\cdot)$ -- \cite[p.471]{GWA}.

We should note that for $r$ of the form $\theta r_o^a$ with $\theta \in \mathbb F_q^\times$, $a \in \mathbb N$
one has 
\[
Se(r,\chi_1,\chi_2,i) = \chi_1(\theta)^i Se(r_o,\chi_1^a,\chi_2,i).
\] 
In particular we can reduce the calculation of $Se(r,\chi_1,\chi_2,i)$ to the case where $r$ is monic. 

\section{The Anderson-Evans-v.Wamelen evaluation of Selberg sums}

In this section we shall summarize the results of Anderson, Evans and v. Wamelen in a notation
suitable for the purposes in hand.   We shall take $r(x)=x^{e_0}(x-1)^{e_1}$.   The group $X_q$ 
is cyclic and the notation of these three authors is based on the choice of a generator denoted
by $\tau$.   The parameters $a$,$b$ and $c$ are related to the parameters here by 
$\tau^a = \chi_1^{e_0}$, $\tau^b = \chi_1^{e_1}$ and $\tau^c = \chi_2$.   It is convenient to 
distinguish two cases.   If $\tau^a = \chi_1^{e_0}$ and $\tau^b = \chi_1^{e_1}$ lie in the subgroup 
generated by $\tau^c = \chi_2$ then we can define two integers $f_0$ and $f_1$ by 
$\chi_1^{e_0}\chi_2^{f_0}=1$ and $\chi_1^{e_1}\chi_2^{f_1}=1$.   It turns out that it is this case, 
the one investigated by v. Wamelen in \cite{vW}, is precisely the one covered by the theory of 
metaplectic forms as formulated in \cite{P2}.   We shall therefore refer to this as the 
\emph{metaplectic} case.   The alternative, when either $\chi_1^{e_0}$ or $\chi_1^{e_1}$ is not 
in the subgroup of $X_q$ generated by $\chi_2$ we call \emph{non--metaplectic}.   This was the 
one investigated by Evans in \cite{RJE2}.   This nomentclature is by no means satisfactory, nor 
even accurate, but it has developed from my private usage and I have not found anything better
to propose.  The distinction does seem to be a useful one as we shall see below.

In the discussion below we shall write $i$ instead of $n$ as used by Anderson, Evans and v.Wamelen.
This corresponds to the usage in Section 2.   We shall also use $n$, the order of $\chi_2$, where 
they used $d$.   For $m \in \mathbb N$ and $j \in \mathbb Z$ we write $(j)_m$ for the least 
non--negative residue of $j$ modulo $m$.

In the theory of Selberg sums a particularly important role is played by a product denoted by 
previous authors by $P_n(a,b,c)$ and which we shall replace by
\[
P_i(e_0,e_1,\chi_1,\chi_2)=\prod_{0 \le j <i} \frac{\tau(\chi_1^{e_0}\chi_2^j)\tau(\chi_1^{e_1}\chi_2^j)\tau(\chi_2^{j+1})\overline{\tau(\chi_1^{e_0+e_1}\chi_2^{i-1+j})}}
{q\tau(\chi_2)}.
\]
This product corresponds to Selberg's product of gamma functions.   It is clear that for a 
certain $A$ one has 
\[
P_{i+\ell n}(e_0,e_1,\chi_1,\chi_2)=A^\ell P_i(e_0,e_1,\chi_1,\chi_2).
\]
One can express $A$ in a simple form, \cite[Lemma 3.1]{RJE2},
\[
\begin{array}{rll}
A&=q^{n-2}\chi_2(\eta(n))/\tau(\chi_2)^n&\mbox{     in the metaplectic case}\\
&=-J(\chi_1^{e_0n},\chi_1^{e_1n})q^{n-1}\chi_2(\eta(n))/\tau(\chi_2)^n&\mbox{    in the non-metaplectic case.}
\end{array}
\]
Note that in the latter case at least one of $\chi_1^{e_0n}$ and $\chi_1^{e_1n}$ is non--trivial.   Our 
notations are the reason that this looks a little different from the formulation in \cite{RJE2}.

It is worth-while noting that in \cite[Lemma 2]{vW} v. Wamelen proves the additional evaluation that
$P_{f_0+f_1+1}(e_0,e_1,\chi_1,\chi_2)$ is equal to
\[
-q^{(f_0+f_1)_n-1}\chi_1(-1)^{f_0+f_1+1}\chi_2(\eta(f_0+f_1+1))
\frac{\tau(\chi_2^{f_0+f_1+1})}{\tau(\chi_2)^{f_0+f_1+1}}.
\]
These evaluations are consequences of the Davenport-Hasse analogue of the Gauss-Legendre multiplication
formula for the gamma function.

For $y \in \mathbb N$ let
\[
T(y,q)= -y + \sum_{k=0}^y (2k+1)(-q)^{y-k}
\]
and
\[
S(y,q)= 1-(1-q)y.
\]
From now on we shall assume the $n>1$ as Anderson's argument has to be modified in the exceptional
case $n=1$.   For more details see \S 5.   The main result of the three authors is that in the non--metaplectic 
case \[Se(x^{e_0}(x-1)^{e_1},\chi_1,\chi_2,i)\] is equal to
\[
\chi_1(-1)^{e_1i}(-1)^i P_i(e_0,e_1,\chi_1,\chi_2)
\]
and in the metaplectic case it is equal to
\[
\chi_1(-1)^{e_1i}(-1)^i q^{\left[ \frac in \right]} P_i(e_0,e_1,\chi_1,\chi_2)
\]
times
\[
\begin{array}{ll} T(2\left[\frac in \right],q) &\mbox{  if } (i)_n \le \mathrm{Min}(f_0,f_1)\le 
\mathrm{Max}(f_0,f_1) < (f_o+f_1-i+1)_n \\
T(2\left[\frac in \right]+1,q) &\mbox{  if } (f_o+f_1-i+1)_n\le \mathrm{Min}(f_0,f_1)\le 
\mathrm{Max}(f_0,f_1) < (i)_n  \\
S(\left[\frac in \right]+1,q) &\mbox{  otherwise.} 
\end{array}
\]

We note here, following \cite[\S 4]{vW}, the following summations,
\[
\sum_{m \ge 0} T(2m,q)X^m = \frac{U_e(q,X)}{(1-X)^2(1-q^2X)}
\]
where
\[
U_e(q,X)=2X^2q^2-qX^2-3qX+X+1,
\]
\[
\sum_{m \ge 0} T(2m+1,q)X^m = \frac{U_o(q,X)}{(1-X)^2(1-q^2X)}
\]
where
\[
U_o(q,X)=X^2q^2+q^2X-3qX-q+2,
\]
and 
\[
\sum_{m \ge 0} S(m,q)X^m = \frac{1+(q-2)X}{(1-X)^2}.
\]
We note also that $U_e(q,q^{-2})= (1-q^{-1})^3$, $U_e(q,1)= 2(q-1)^2$, $U_o(q,q^{-2})= -q(1-q^{-1})^3$ and $U_o(q,1)= 2(q-1)^2$.
These results mean that in the non-metaplectic case and in any complex embedding we have
\[
\mltsum{i\ge 0}{i\equiv i_0} S(x^{e_0}(x-1)^{e_1},\chi_1,\chi_2,i) X^{(i-i_0)/n}
\]
is equal to
\[
\chi_1(-1)^{e_1i_0}(-1)^{i_0}P_{i_0}(e_0,e_1,\chi_1,\chi_2)/(1-\chi_1(-1)^{e_1n}(-1)^nAX)
\]
where we assume that $0\le i_0<n$ and $A$ is as above.   Note that $|A|$ takes on, if $n\neq 1$,
one of the three values $q^{\frac n2 -1}$, $q^{\frac n2 -\frac 12}$ or $q^{\frac n2}$ in any 
complex embedding.   If $n=1$ then $|A|$ takes on one of the three values $1$, $q^{\frac 12}$ or 
$q$ in any complex embedding.

In the metaplectic case we find that, according to the three cases above, the series is equal to 
\[
\chi_1(-1)^{e_1i_0}(-1)^{i_0}\frac{P_{i_0}(e_0,e_1,\chi_1,\chi_2)U_e(q,\chi_1(-1)^{e_1n}(-1)^nAX)}
{(1-\chi_1(-1)^{e_1n}(-1)^n qAX)^2(1-\chi_1(-1)^{e_1n}(-1)^nq^3AX)},
\]
\[
\chi_1(-1)^{e_1i_0}(-1)^{i_0}\frac{P_{i_0}(e_0,e_1,\chi_1,\chi_2)U_0(q,\chi_1(-1)^{e_1n}(-1)^n)AX}
{(1-\chi_1(-1)^{e_1n}(-1)^n qAX)^2(1-\chi_1(-1)^{e_1n}(-1)^nq^3AX)},
\]
or 
\[
\chi_1(-1)^{e_1i_0}(-1)^{i_0}\frac{P_{i_0}(e_0,e_1,\chi_1,\chi_2)(1+(q-2)\chi_1(-1)^{e_1n}(-1)^n qAX)^2}
{(1-\chi_1(-1)^{e_1n}(-1)^n qAX)^2}
\]
respectively.  Note that the nature of the singularities reflect the type of the Selberg sum.   As we have assumed that 
$n>1$ we have in the metaplectic case $|A|=q^{\frac n2 - 2}$.

\section{The transformation formula}
We shall now turn to a property of Selberg sums analogous to Theorem 1 of \cite{P2}.   
We need some preparations in order to be able to formulate the result.   Let $n$ be 
the order of $\chi_2$ and let $n'$ be the order of $\chi_1$.   Let $\chi_0$ be such 
that $\chi_1$ and $\chi_2^2$ are in the group generated by $\chi_0$.   Let $a\ge 0$ and 
$b\ge 0$ be such that
\[
\chi_1=\chi_0^a, \mbox{  } \chi_2^2=\chi_0^b.
\]

Let $\pi$ monic and a prime in $R=\mathbb F_q[x]$.   If $\pi | r$ then
\[
Se(\pi^{n'}r,\chi_1,\chi_2,i) = Se(r,\chi_1,\chi_2,i).
\]
If $\pi \not| r$ this is no longer true.   We have that 
\[
Se(r,\chi_1,\chi_2,i)-Se(\pi^{n'}r,\chi_1,\chi_2,i)
\]
is equal to 
\[
\sum_{\mathrm{deg}(c_1\pi)=i} \mu(c_1\pi)\chi_1(r/c_1\pi)\chi_2(D(c_1\pi)).
\]
Since $D(c_1\pi)=D(c_1)D(\pi)R(c_1,\pi)^2$ this becomes
\[
-\chi_1(r/\pi)\chi_2(D(\pi))\sum_{\mathrm{deg}(c_1)=i-\mathrm{deg}(\pi)} \mu(c_1)\chi_1(r/c_1\pi)\chi_2(\pi/c_2)^2\chi_2(D(c_1))
\]
or 
\[
-\chi_1(r/\pi)\chi_2(D(\pi)Se(r^a\pi^{b},\chi_0,\chi_2,i-\mathrm{deg}(\pi)).
\]
We regard this as a stability property of Selberg sums.

\noindent\textbf{Theorem 1}
\emph{Let $\alpha,\beta,\gamma,\delta \in \mathbb F_q$ be such that $\Delta =\alpha\delta - -\beta\gamma \neq 0$.
Then for suitable integers $M,M'$ we have that} 
\[
\chi_2(\Delta)^{1-i}Se\left(r\left(\frac{\alpha x + \beta}{\gamma x + \delta}\right)^a\left(\frac 1{(\gamma x + \delta)^2}
\right)^{b(i-1)+M},\chi_0,\chi_2,i\right) 
\]
\emph{is equal to}
\[Se(r(x)(-\gamma x + \alpha)^{M'},\chi_1,\chi_2,i)
\]
\emph{If $\gamma=0$ we can take $M,M'=0$.   If $\gamma \neq 0$ then $M$ is to such that $\chi_0^M=1$
and $M+b(i-1) > a\cdot\mathrm{deg}(r)$, and $M'$ is such that $M'>0$ and $\chi_1^{M'}=1$.} 

\noindent\emph{Proof:} The proof is carried out by verifying the identity for $\left(\begin{array}{cc}\alpha&\beta\\
\gamma&\delta \end{array}\right)$ of the form $\left(\begin{array}{cc}1&\lambda\\ 0&1 \end{array}\right)$,
of the form $\left(\begin{array}{cc}\alpha&0\\0 &\delta \end{array}\right)$ and equal to 
$\left(\begin{array}{cc}0&1\\ 1&0 \end{array}\right)$ and then combining these in the usual manner.

The first two cases are straightforward; we need only replace $c$ in the sum defining the Selberg
sum by $c(x+\lambda)$ and $c(\theta x)\theta^{-i}$ where $\theta = \alpha/\delta$ respectively.

In the case  $\left(\begin{array}{cc}\alpha&\beta\\ \gamma&\delta \end{array}\right)=
\left(\begin{array}{cc}0&1\\ 1&0 \end{array}\right)$ we assume that $r(0)=0$; this is the
reason for the introduction of the parameter $M'$.   This being so $c$ is not divisible by $x$
and so $c(0)\neq 0$.   We replace $c$ in the sum by $c(x^{-1})x^i/c(0)$ and recall that 
$R(x,c)=c(0)$.   The formula quoted now follows from the theory of resultants. 

\section{Selberg sums and metaplectic groups}
The results of Section 3 show that there is a large number of relationships between Selberg
sums.    The results of Anderson, Evans and v.~Wamelen show that polynomials with two rational
zeros over $\mathbb F_q$ can be evaluated explicitly.   It turns out, as we shall see in the 
next section, that if the zeros are no longer rational a similar formula holds.    The 
next case which one can investigate, as in \cite{P2} is that of polynomials with three 
rational zeros.   In view of Theorem 1 these can be brought to the form 
$x^{e_0}(x-1)^{e_1}(x-\lambda)^{e_\lambda}$.   As in \cite{P2} there are a large number 
of relations between these, now a little more complicated as we have to move into the
region of stability.   They are the analogues of the relations beween such sums 
similar to transformations of hypergoemetric functions (cf. \cite[Chap. XIV]{W&W}).   We should 
note that in this standard form such Selberg sums can be easily computed as the resultant
can be evaluated by a continued-fraction type of recursion.   This is implement in,
for example, the gp/PARI package.   I hope to discuss the results of these calculations
in a future paper.

There are some further relations that are of interest.   Let $n,n'$ be as in Section 3.
If $n'|n$ then the generating series
\[
\mltsum{i\equiv i_0 \pmod n}{i\ge 0} Se(r,\chi_1,\chi_2,i)T^{(i-i_0)/n}
\] 
which converges in $|T|<q^{-n}$, is, for each $i_0: 0\le i <n$, a rational function in $T$.    There 
is at most one singularity in $|T|<q^{-n/2}$ and this, if it exists, is simple and is located at 
$T=q^{-1}(-\tau(\chi_2^{-1}))^{-n}$.   These are the series investigated in \cite{P1} and \cite{P2}.  
The determination of the residue and the establishing of its properties is one of the major questions 
addressed in those papers and it is one that has, as yet, only been partially answered.

If the condition $n'|n$ is not satisifed then the series above can also be investigated by
means of the theory of Eisenstein series.   The point here is that one has to use an
Eisenstein series of a ``Nebentypus'' depending on $\chi_1$ and $\chi_2$.   The general 
theory in contained in \cite{KP}.   In fact in this case the constant terms of the Eisenstein
series will be made up from holomorphic L-series.   The analytic continuation then follows 
from the ``principle of the constant term''; see \cite[Theorem 4.8]{LM} for a statement and 
\cite[Theorem 1.6.6]{GH} for a proof in the function--field case.   The details have not 
been given explicitly but it seems very plausible that in this case one will be able to 
conclude that the series above  also in this case represents a rational function but one
now with no singularities in $|T| < q^{-n/2}$.   One expects that one will be able to 
determine the denominator and to estimate the degree of the numerator which are known in the 
metaplectic case also in the non-metaplectic case.

In is instructive to examine the special case $n=1$.   There are two cases, the metaplectic case 
in which $\chi_1(r/\cdot)$ is a principal character and the other case when $\chi_1(r/\cdot)$ is 
non--principal.   In the first case we see that $ \sum_{i\ge 0} Se(r,\chi_1,\chi_2,i)T^{i}$ 
is $L(s,\chi_1(r/\cdot)^{-1}$ with $T=q^{-s}$.   Let $r_0$ be the modulus of $\chi_1(r/\cdot)$;
then this series is equal to
\[
(1-q^{1-s})\prod_{\pi | r_0}(1-q^{-\mathrm{deg}(\pi )s}).
\] 
In this case the generating series is a polynomial.

If $\chi_1(r/\cdot)$ is non--principal then, if we suppose for the moment that the $\chi_1(r/\cdot)$ is 
primitive.   Then Anderson, \cite{GWA}, has observed that $L(s,\chi_1(r/\cdot)$ is a polynomial 
in $q^{-s}$ of degree $\mathrm{deg}(r_0)-1$.   It follows that $L(s,\chi_1(r/\cdot)^{-1}$ has this 
number of singularities, counted with multiplicities.   They lie on $q^{-s}=q^{-1/2}$ by Weil's theorem.
If the character is not primitive then there is a numerator of the form 
$\prod(1-q^{-s\mathrm{deg}(\pi)})$ where the product is over all the primes dividing the 
modulus but not the conductor.   In the case of the standard Selberg sums the degree of is
$1$ and it is of the form $1-\alpha q^{-s}$ where $\alpha$ is a Jacobi sum.  The case
of other $r$ will be much less straightforward.  
  
We return to the general case.   The singularities of the generating series of the sequence 
$Se(r,\chi_1,\chi_2,i)$ are given as certain sums of their coefficients, i.e. of Selberg sums.
This has been described in \cite{P1},\cite{P2}.   The theory of metaplectic groups then leads 
to a number of relations between the residues, the so--called ``Hecke relations''.   A 
consequence of this is that there are a number of unexpected relations between Selberg 
sums.   Some examples are given in \cite[\S 4]{P1}.   To date there is no systematic
method of treating these relationships.    This is not unlike the situation in the 
theory of generalized hypergeometric functions.

It is interesting to consider Anderson's method in this connection.   His approach is 
based on a formula which is derived from the theory of Dirichlet series over $R$.   We
shall now sketch his technique.   Let $f$ be a monic polynomial in $R$.   Let $\chi$ be
as before and let $n$ be its order.   We regard $\chi(f/\cdot)$ as a Dirichlet character.   
Let $f_o$ be the conductor of $f$.   This is $\prod \pi$ where $\pi$ runs through those
monic prime divisors of $f$ whose order in $f$ is not divisible by $n$.  We shall 
assume that $f_o \neq 1$.   Then by \cite[Prop. 2.1]{GWA}
\[
\omega(D(f_o))\chi(f/f_o')\tau(\chi^{-1})^{\mathrm{deg}(f)-1}
\prod_{\pi | f_o} \tau(\chi^{-\mathrm{ord}_\pi(f)})^{-\mathrm{deg}(\pi)}
\]
is equal to 
\[
\mltsum{g: \mathrm{deg}(g)=\mathrm{deg}(f_o)-1}{g \mbox{  \scriptsize{monic}}} \chi(f/g).
\]
This is proved by using the functional equation for the L-function 
$L(s,\chi(f/\cdot))$ and the fact that the latter is a polynomial to determine the 
coefficient of $(q^{-s})^{\mathrm{deg}(f_o)-1}$.   Anderson proves this in the case $n=q-1$
but the proof is valid whenever $\chi(f/\cdot)$ is not principal.   Anderson achieves the 
same generality through his parameter $c$.

Anderson's proof of the formula for the Selberg sum exploits the evaluation of a double sum 
in two different ways.   The argument is strongly reminiscient of the the multiple Dirichlet
series technique, as used, for example, in \cite{CM}.   The crucial point of the argument is that 
$r_o$ should be quadratic and it follows that one can also evaluate Selberg sums explicitly 
when $r$ is a power of an irreducible quadratic polynomial.

It seems plausible that one could use Anderson's method combined with elementary considerations
to determine the properties of metaplectic Eisenstein series over rtional function fields.   In
fact the ``principle of the constant term'' referred to above is relatively elementary and as 
the constant term is easy to study in this case even the Eisenstein series approach is relatively
elementary.   The method of multiple Dirichlet series is, at least for higher ranks, an important
component in the study of Eisenstein series and so the applicibility of both methods should not be 
surprising.   At any rate it should reassure those unfamiliar with this theory that the method is, 
at heart, in this case at least, elementary.

\noindent
Mathematisches Institut\\
Bunsenstr. 3--5\\
37073 G\"ottingen\\\
Germany

\noindent
e-mail:\texttt{spatter@gwdg.de}


\begin{thebibliography}{99}

\bibitem{GWA} G.W.Anderson: The evaluation of Selberg 
sums,\textit{Comptes Rendus Acad.Sci.Paris,Ser. I}, 311(1990)469-472.

\bibitem{CR} D. Carmon, Z.Rudnick, The autocorrelation function of the M\"obius function
and Chowla's conjecture for the ration function field, \textit{Quart. J. Math.}65(2014)53--61.

\bibitem{CM}G. Chinta, J.B. Mohler, Sums of L-functions over rational function fields,
\textit{Acta Arith.} 144(210)53--68. 

\bibitem{KC} K.Conrad: Irreducible values of polynomials: a non-analogy,in
\emph{Number fields and function fields -- two parallel worlds} Birkh\"auser,2005,
71--86.

\bibitem{RJE1} R.J. Evans: Identities for products of Gauss sums over finite fields, 
\textit{L'Enseign.Math.} 27(1981)197-209.

\bibitem{RJE2} R.J. Evans: The evaluation of Selberg character sums, 
\textit{L'Enseign.Math.} 37(1991)235-248.

\bibitem{GH} G.Harder: Chevalley groups over function fields and automorphic forms,
\textit{Ann.Math.} 100 (1974) 249-306.

\bibitem{KP} D.A. Kazhdan, S.J. Patterson: Metaplectic Forms,
\textit{Publ.Math. IHES} 59 (1984) 35-142.

\bibitem{LM} L.E.Morris, Eisenstein series for reductive groups 
over global function fields, The Cusp Form Case, 
\emph{Can. J. Math.}34(1982)91-168.

\bibitem{P1} S.J. Patterson, Note on a paper of J. Hoffstein,
\textit{Glasgow Math. J.} 296 (2007) 125-161, 217-220.

\bibitem{P2} S.J. Patterson, The Fourier coefficients of metaplectic 
theta series on GL(2) over rational function fields, \emph{in preparation}

\bibitem{ASCW1} A. Selberg,\textit{ Collected Papers, 1}, Springer, 1989.

\bibitem{RGS} R. G. Swan, Factorization of polynomials over finite fields, 
\textit{Pacific J. Math.} 12(1962)1099--1106.

\bibitem{vdW} B.L. van der Waerden: \textit{Modern Algebra, Vol. 1},
(Trans. F. Blum) Fredrick Unger  Publishing Co.,New York, 1949.

\bibitem{vW} P.B.van Wamelen: Proof of the Evans-Root conjectures for 
Selberg character sums, \textit{J.Lond.Math.Soc.,II.Ser.}48(1993)415-426.

\bibitem{BNT} A. Weil, \textit{Basic Number Theory}, Springer, 1967.

\bibitem{W&W} E.T.Whittaker, G.N.Watson, \textit{Modern Analysis}, 4th edition, Cambridge 
Univ. Press, 1946
\end{thebibliography}
\end{document}